\documentclass[10pt,twoside]{article}
\usepackage{amsfonts}
\usepackage{amssymb}
\usepackage{epsfig}
\usepackage{qtree}
\input xypic
\xyoption{all}
\LaTeXdiagrams
\xyoption{2cell}
\UseAllTwocells
\newtheorem{defi}{Definition}
\newtheorem{teo}{Theorem}

\def\be{\begin{equation}}
\def\ee{\end{equation}}
\def\raa{\rightarrow}
\def\lra{\longrightarrow}

\def\circt{\ddot{\circ}}
\def\N{{\mathbb N}}

\def\b{\fbox{\bf B}}
\def\c{\fbox{\bf C}}
\def\r{\fbox{\bf R}}

\def\PRdesc{{\bf PRdesc}}
\def\PRalg{{\bf PRalg}}
\def\PRfunc{{\bf PRfunc}}
\def\PRalgC{{\bf PRalgC}}
\def\PRalgI{{\bf PRalgI}}
\def\PRalgCat{{\bf PRalgCat}}
\def\PRalgCatX{{\bf PRalgCatX}}
\def\PRalgCatN{{\bf PRalgCatN}}
\def\PRalgCatXN{{\bf PRalgCatXN}}

\def\la{\langle}
\def\ra{\rangle}
\newenvironment{examp}{ \stepcounter{examnum} {\bf \noindent Example:}}{$\Box$}

\newcounter{examnum}[section]
\newcounter{remarnum}[section]
\begin{document}
\title{Galois Theory of Algorithms}
\author{Noson S. Yanofsky\footnote{Department of Computer and Information Science,
Brooklyn College CUNY,
Brooklyn, N.Y. 11210. And Computer Science Department,
The Graduate Center CUNY,
New York, N.Y. 10016.
e-mail: noson@sci.brooklyn.cuny.edu}}
\maketitle

\centerline{\it  In honor of Rohit Parikh}

\begin{abstract}
\noindent Many different programs are the implementation of the same
algorithm. The collection of programs can be partitioned into different classes corresponding to the algorithms they implement. This makes the collection of algorithms a quotient of the collection of programs. Similarly, there are many different algorithms
that implement the same computable function. The collection of algorithms can be partitioned into different classes corresponding to what computable function they implement. This makes the collection of computable functions into a quotient of the 
collection of algorithms.
Algorithms are intermediate between programs and functions:

    Programs $\twoheadrightarrow$ Algorithms $\twoheadrightarrow$ Functions.

\noindent Galois theory investigates the way that a subobject sits inside an object. We investigate how a quotient object sits inside an object. By looking at the Galois group of programs, we study the intermediate types of algorithms possible and the types of structures these algorithms can have.
\end{abstract}

\section{Introduction}

As an undergraduate at Brooklyn College in 1989, I had the good fortune to take a masters level course in theoretical computer science given by Prof. Rohit Parikh. His infectious enthusiasm and his extreme clarity turned me onto the subject. I have spent the last 25 years studying theoretical computer science with Prof. Parikh at my side. After all these years I am still amazed by how much knowledge and wisdom he has at his fingertips. His broad interests and warm encouraging way has been tremendously helpful in many ways. I am happy to call myself his student, his colleague, and his friend. I am forever grateful to him. 

In this paper we continue the work in \cite{Yano} where we began a study of formal definitions of algorithms
(knowledge of that paper is not necessary for this paper.) The previous paper generated some interest
in the community: Yuri I. Manin looked at the structure of programs and algorithms from the operad/PROP point
of view \cite{Manin} (Chapter 9), see also \cite{Manin1, Manin2} where it is discussed in the context of renormalization;
there is an ongoing project to extend this work from
primitive recursive functions to all recursive functions in \cite{ManinYano}; Ximo Diaz Boils has looked
at these constructions in relations to earlier papers such as \cite{Burroni,Maietti,mft}; Andreas Blass, Nachum Dershowitz, and Yuri Gurevich discuss the paper in \cite{BlassDerGur} with reference to their definition of an algorithm.

\begin{figure}[!h]
\centering
  \includegraphics[height=3.25 in]{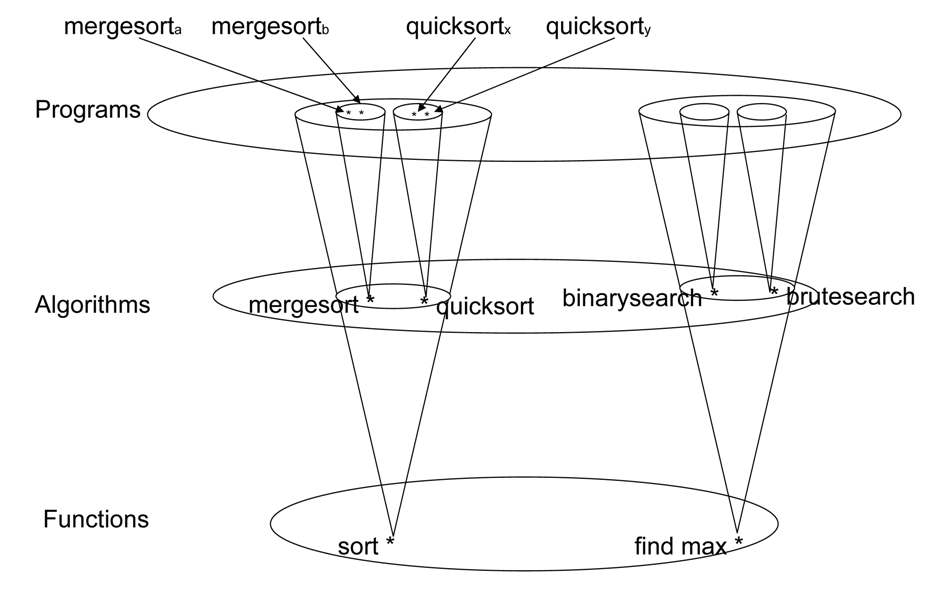}\\
\caption{Programs, Algorithms and Functions.}
\label{firstpic}
\end{figure}
Figure \ref{firstpic} motivates the
formal definition of algorithms. On the bottom is the set of computable functions. Two examples of computable functions are given:
the sorting function and the find max function. On top of the diagram is the set of
programs. To each computable function on the bottom, a cone shows the corresponding set of
programs that implement that function. Four such programs that implement the sorting function have been highlighted:
$\verb"mergesort"_a$, $\verb"mergesort"_b$, $\verb"quicksort"_x$ and $\verb"quicksort"_y$. One can think of
$\verb"mergesort"_a$ and $\verb"mergesort"_b$ as two implementations of the mergesort algorithm
written by Ann and Bob respectively. These two programs are obviously similar but they are not the same.
In the same way, $\verb"quicksort"_x$ and $\verb"quicksort"_y$ are two different implementations of the quicksort
algorithm. These two programs are similar but not the same. We shall discuss in what sense they are ``similar.'' Nevertheless programs that
implement the mergesort algorithm are different than programs that implement the quicksort algorithm.
This leads us to having algorithms as the middle level of Figure \ref{firstpic}. An algorithm is to
be thought of as an equivalence class of programs that implement the same function. The mergesort algorithm
is the set of all programs that implement mergesort. Similarly, the quicksort algorithm is the set of all
programs that implement quicksort.
The set of all programs are partitioned into equivalence classes and each equivalence class corresponds to an
algorithm. This gives a surjective map from the set of programs to the set of algorithms.

One can similarly partition the set of algorithms into equivalence classes. Two algorithms
are deemed equivalent if they perform the same computable function. This gives a surjective function
from the set of algorithms to the set of computable functions.

This paper is employing the fact that equivalence classes of programs
have more manageable structure than the original set of programs. We will find that the set of programs does not have much structure at all.
In contrast, types of algorithms have better structure and the set of computable functions have
a very strict structure.

The obvious question is, what are the equivalence relation that say when two programs are ``similar?'' In \cite{Yano}
a single tentative answer was given to this question. Certain relations were described that seem universally agreeable. Using
these equivalence relations, the set of algorithms have the structure of a category (composition)  with a product (bracket) and
a natural number object (a categorical way of describing recursion.) Furthermore, we showed that with these
equivalence relations, the set of algorithms has universal properties. See \cite{Yano} for more details.

Some of the  relations that describe when two programs are ``similar'' were:
\begin{itemize}
\item One program might perform $Process_1$ first and then perform an unrelated $Process_2$ after.
The other program might perform the two unrelated processes in the opposite order.
\item One program will perform a certain process in a
loop $n$ times and the other program will ``unwind the loop'' and perform it $n-1$ times and then perform the process again
outside the loop.
\item One program might perform two unrelated processes in one loop, and the other
program might perform each of these two processes in its own loops.
\end{itemize}

In \cite{BlassDerGur}, the subjectivity of the question as to when two programs
are considered equivalent was criticized. While writing \cite{Yano}, we were aware that the answer to this question is a subjective decision
(hence the word ``Towards'' in the title), we nevertheless described the structure of algorithms in that particular
 case. In this paper we answer that assessment of \cite{Yano} by looking at the {\it many different sets} of equivalence relations that
one can have. It is shown that with every set of equivalence relations we get a certain structure.

The main point of this paper is to explore the set of possible intermediate structures between programs and computable functions using the
techniques of Galois theory. In Galois theory, intermediate fields are studied by looking at automorphism of fields.
Here we study intermediate algorithmic structures by looking at automorphism of programs.

A short one paragraph review of Galois theory is in order. Given a polynomial
with coefficients in a field $F$, we can ask if there is a solution to the polynomial in an extension field $E$.
One examines the group of automorphisms of $E$ that fix $F$, i.e., automorphisms $\phi: E \lra E$ such that for all
$f\in F$ we have $\phi(f)=f$.
\be \xymatrix{
E\ar@{>>}[ddddrrrr] \ar[rrrrrrrr]^\phi &&&&&&&&E
\\
\\
&&& F_H \ar[rr]^\psi \ar@{^{(}->}[uulll]\ar@{^{(}->}[uurrrrr]&& F_H
\\
\\
&&&&F\ar@{^{(}->}[uuuullll] \ar@{^{(}->}[uuuurrrr]\ar@{^{(}->}[uul]\ar@{^{(}->}[uur]
\label{classicalGal}}\ee
This group is denoted $Aut(E/F)$. For every subgroup $H$ of
$Aut(E/F)$ there is an intermediate field $F\subseteq F_H \subseteq E$. And conversely for every
intermediate field $F\subseteq K \subseteq E$ there is a subgroup $H_K=Aut(E/K)$ of $Aut(E/F)$. These two maps form a Galois connection between intermediate fields and subgroups of  $Aut(E/F)$. If one further restricts to normal intermediate fields and normal subgroups than there is an isomorphism of partial orders.  
This correspondence is the essence of the Fundamental Theorem of Galois Theory which says that the lattice of normal subgroups of $Aut(E/F)$
is isomorphic to the dual lattice of normal intermediate fields between $F$ and $E$. The properties of $Aut(E/F)$ mimic the
properties of the fields. The group is ``solvable'' if and only if the polynomial is ``solvable.''

In order to understand intermediate algorithmic structures
we study automorphisms of programs. Consider all automorphisms of programs that
respect functionality. Such automorphisms can be thought of as ways of swapping programs for other programs that perform the same function. An automorphism $\phi$ makes the following outer triangle commute.
\be\xymatrix{ \label{GalofAlg}
\mbox{Prgs}\ar@{>>}[ddddrrrr]_Q \ar@{>>}[ddrrr]^R\ar[rrrrrrrr]^\phi &&&&&&&&\mbox{Prgs}\ar@{>>}[ddddllll]^Q\ar@{>>}[ddlllll]^R
\\
\\
&&& \mbox{Algs} \ar@{>>}[ddr]^S \ar[rr]^\psi && \mbox{Algs} \ar@{>>}[ddl]_S
\\
\\
&&&&\mbox{Fncts}
}\ee
A subgroup of the group of all automorphisms is going to correspond to
an intermediate structure. And conversely, an intermediate algorithmic structure will correspond to a subgroup. We will then consider special types of ``normal'' structures to get an isomorphism of partial orders. This will be the
essence of the fundamental theorem of Galois theory of algorithms. This theorem formalizes the intuitive notion that two programs
can be switched for one another if they are considered to implement the same algorithm. One extreme case is if you consider every program to be its own algorithm. In that case
there is no swapping different programs. The other extreme case is if you consider two programs to be equivalent when they perform the same
function. In that case you can swap many programs for other programs. We study all the intermediate possibilities.

Notice that the diagonal arrows in Diagram (\ref{GalofAlg}) go in the opposite direction of the arrows in Diagram (\ref{classicalGal})
and are surjections rather
than injections. In fact the proofs in this paper are similar to the ones in classical Galois theory as long as you stand on your head.
We resist the urge to call this work
``co-Galois theory.''

All this is somewhat abstract. What type of programs are we talking about?
What type of algorithmic structures are we dealing with? How will our descriptions be specified?
Rather than choosing one programming language to the exclusion of others, we look at a
language of descriptions of primitive recursive functions. We choose this language because
of its beauty, its simplicity of presentation, and the fact that most readers are familiar with this language.
The language of descriptions of primitive recursive functions has only three operations: Composition, Bracket, and Recursion.
We are limiting ourselves to the set
of primitive recursive functions as opposed to all computable functions for ease. By so doing, we are going to get a proper subset of all algorithms. Even though we are, for the present time, restricting
ourselves, we feel that the results obtained are interesting in their own right. There is an ongoing project to extend this
work to all recursive functions \cite{ManinYano}.

\vspace{.2in}
Another way of looking at this work is from the homotopy theory point of view. We can think of the set of
programs as a graph enriched over groupoids. In detail, the 0-cells are the powers of the natural number (types),
the 1-cells are the programs from a power of natural numbers to a power of natural numbers.
There is a 2-cell from one program to another program if and only if they are ``essentially the same''. That is,
the 2-cells describe the equivalence relations. By the symmetry of the equivalence relations, the 2-cells form a
groupoid. (One goes on to look at the same graph structure with different enrichments. In other words, the 0-cells and the 1-cells are the same, but look at different possible isomorphisms between the 1-cells.) Now we take the quotient, or fraction category where we identify the programs at the end
of the equivalences. This is the graph or category of algorithms. From this perspective we can promote the mantra:

\vspace{.2in}

\centerline{``Algorithms are the homotopy category of programs.''}

\vspace{.2in}
Similar constructions lead us to the fact that 

\vspace{.2in}

\centerline{``Computable functions are the homotopy category of algorithms.''}

\vspace{.2in}
\noindent This is a step towards
\vspace{.2in}

\centerline{``Semantics is the homotopy category of syntax.''}

\vspace{.2in}
\noindent Much work remains to be done.

\vspace{.2in}

Another
way of viewing this work is about composition.
Compositionality has for many decades been recognized as one of the
most valuable tools of software engineering.

There are different levels of abstractions that we use when we teach
computation or work in building computers,
networks, and search engines. There are programs, algorithms, and functions.
Not all levels of
abstraction of computation admit useful structure.  If we take
programs to be the finest level then we may find it hard to compose
programs suitably.  But if we then pass to the abstract functions they
compute, again we run into trouble. In between these two extremes ---extreme
concreteness and extreme abstractness--- there can be many levels of
abstraction that admit useful composition operations unavailable at
either extreme.

It is our goal here to study the many different levels of algorithms and
to understand the concomitant different possibilities of composition.
We feel that this work can have great potential value for software engineering.

\vspace{.2in}

Yet another way of viewing this work is an application and a variation of some ideas from universal algebra and model theory.
In the literature, there is some discussion of Galois theory for arbitrary universal algebraic structures (\cite{Cohen} section II.6) and and m-theoretic
structures  (\cite{daCostaRodrigues,daCosta,Fleisher,Poschel}.)
In broad philosophical terms, following the work of Galois and Klein's {\it erlangen} program, an object can be defined
by looking at its symmetries. Primitive recursive programs are here considered as a universal algebraic structure where the generators of the
structure are the initial functions while composition, bracket and recursion are the operations.
This work examines the symmetries of such programs and types of structures that can be
defined from those symmetries.

\vspace{.2in}

Section 2 reviews  primitive recursive programs and the basic structure that they have. In Section 3 we define
an algorithmic universe as the minimal structure that a set of algorithms can have. Many examples are given. The main theorems in this paper
are found in Section 4 where we prove the Fundamental Theorem of Galois Theory for Algorithms. Section 5 looks at our work from the point of view of homotopy theory, covering spaces, and Grothendieck Galois theory. We conclude with a list of possible ways that this work might progress in the future.

\vspace{.2in}

\noindent {\bf Acknowledgment}. I thank Ximo Diaz Boils, Leon Ehrenpreis (of blessed memory), Thomas Holder, Roman Kossak, Florian Lengyel, Dustin Mulcahey, Robert Par\'{e}, 
Vaughan Pratt, Phil Scott, and Lou Thrall for helpful discussions. I am also thankful to an anonymous reviewer who was very helpful. 

\section{Programs}

Consider the structure of all descriptions of primitive recursive functions. Throughout this paper we use the words
``description'' and ``program'' interchangeably. The descriptions form a graph
denoted $\PRdesc$. The objects (nodes) of the graph are powers of natural numbers $\N^0, \N^1, \N^2, \ldots, \N^n, \ldots$
and the morphisms (edges) are descriptions of primitive recursive functions. In particular, there exists descriptions of initial
functions: the null function (the function that only outputs a 0) $z:\N^0={\bf 1} \lra \N$, the successor function $s:\N \lra \N$, and the projection functions, i.e.,
for all $n \in \N$ and for all $1\le i \le n$ there are distinguished descriptions $\pi^n_i:\N^n \lra \N$.

There will be three ways of composing edges in this graph:
\begin{itemize}
\item Composition: For $f:\N^a \lra \N^b$ and $g:\N^b \lra \N^c$, there is a $$(g \circ f):\N^a \lra \N^c.$$
Notice that this composition need not be associative. There is also no reason to assume that this composition has a unit.
\item Recursion: For $f:\N^a \lra \N^b$ and $g:\N^{a+b+1} \lra \N^b$, there is a $$(f \sharp g):\N^{a+1} \lra \N^b. $$
There is no reason to think that this operation satisfies any universal properties or that it respects the composition or the bracket.
\item Bracket: For $f:\N^a \lra \N^b$ and $g:\N^a \lra \N^c$, there is a $$\la f, g \ra:\N^a \lra \N^{b+c}.$$ There is no
reason to think that this bracket is functorial (that is respects the composition) or is in any way coherent.

\end{itemize}

At times we use trees to specify the descriptions. The leaves of the trees will have initial functions and the internal nodes will be
marked with {\bf C}, {\bf R} or {\bf B} for composition, recursion and bracket as follows:

\qtreecenterfalse
\Tree[{$f:\N^a\raa \N^b$} {$g:\N^b\raa \N^c$} ] .{$g\circ f:\N^a \raa \N^c$\\
\c} \qquad  \Tree[ {$f:\N^a\raa \N^b$} {$g:\N^a \times \N^b \times \N \raa \N^b$} ] .{$f\sharp
g:\N^a\times \N \raa \N^b$\\ \r} \qquad  \Tree[ {$f:\N^a\raa \N^b$} {$g:\N^a\raa
\N^c$} ] .{$\la f, g \ra :\N^a \raa \N^b \times \N^c$\\ \b}

Just to highlight the distinction between programs and functions, it is important to realize that the following are all
legitimate descriptions of the null function:
\begin{itemize}
\item $z:\N^{0} \lra \N$
\item $(z\circ s \circ s \circ s \circ  s \circ z \circ s \circ s \circ s \circ s \circ s \circ s \circ s \circ  s \circ s \circ z): \N^{0} \lra \N$
\item $ (z \circ (\pi^2_1 \circ \la s, s\ra))\circ z:\N^{0} \lra \N$
\item etc.
\end{itemize}
There are, in fact, an infinite number of descriptions of the null function.

\vspace{.2in}
In this paper we will need ``macros'', that is, certain combinations of operations to get commonly used descriptions.  Here are a few.

There is a need to generalize the notion of a projection. The $\pi^n_i$
accepts $n$ inputs and outputs one. A multiple projection takes $n$ inputs and outputs $m$ outputs.
Consider $\N^n$ and the sequence
$X= \la x_1,x_2, \ldots , x_m \ra$ where each $x_i$ is in $\{1,2, \ldots, n \}$. For every $X$ there exists
$\pi^{X}:\N^n \raa \N^m$ as
$$ \pi^X = \la \pi^n_{x_1}, \la \pi^n_{x_2}, \la \ldots, \la \pi^n_{x_{m-1}},\pi^n_{x_m} \ra \ra \ldots \ra.$$
In other words, $\pi^X$ outputs the proper numbers in the order described by $X$.
In particular
\begin{itemize}
\item If $I=\la 1,2,3, \ldots , n\ra$ then $\pi^I:\N^n \lra \N^n$ will be a description of the identity function.
\item If $X=\la 1,2,3, \ldots , n,1,2,3, \ldots, n\ra$ then $\pi^X=\triangle=\pi^n_{n+n} : \N^n \lra \N^{n+n}$ is the diagonal map.
\item For $a\le b \in \N$, if $$X=\la b, b+1, b+2, \ldots, b+a, 1, 2, 3, \ldots, b-1 \ra,$$
then $\pi^X$ will be the twist operator which
swaps the first $a$ elements with the second $b$ elements. Then by abuse of notation, we shall write $$\pi^X=tw=\pi^{a+b}_{b+a}:\N^{a+b} \lra \N^{a+b}.$$
\end{itemize}
Whenever possible, we omit superscripts and subscripts.
\vspace{.2in}

Concomitant with the bracket operation is the product operation. A product of two maps is defined for a given $f:\N^a \raa\N^b$ and $g:\N^c \raa\N^d$
as
$$f \times g:\N^a\times \N^c  \raa\N^b \times \N^d.$$ The product can be defined using the bracket as
$$f \times g = \la f \circ \pi^{a + c}_a, g \circ \pi^{a + c}_c \ra.$$

Given the product and the diagonal, $\triangle=\pi^a_{a+a}$, we can define the bracket as
$$\xymatrix{
\N^a \ar[rr]^{\la f,g \ra} \ar[rdd]_\triangle  && \N^b \times \N^c
\\
\\
& \N^a \times \N^a. \ar[ruu]_{f\times g}
}$$

Since the product and the bracket are derivable from each other, we use them interchangeably.

That is enough about the graph of descriptions. 

\vspace{.2in}

Related to descriptions of primitive recursive functions
is the set of primitive recursive {\it functions}. The set of functions has a lot more structure than
$\PRdesc$. Rather than just being a graph, it forms a category.
$\PRfunc$ is the category of primitive recursive functions. The objects of this category
are powers of natural numbers $\N^0, \N^1, \N^2, \ldots, \N^n, \ldots$
and the morphisms are primitive recursive functions. In particular, there are specific maps
$z:\N^{0} \lra \N$, $s:\N \lra \N$ and for all $n \in \N$ and for all $1\le i \le n$ there
are projection maps $\pi^n_i:\N^n \lra \N$. Since composition of primitive recursive functions is associative
and the identity functions $id=\pi^n_n:\N^n \lra \N^n$ are primitive recursive and act as units for composition, $\PRfunc$ is a
genuine category. $\PRfunc$ has a categorically coherent Cartesian product $\times$.
Furthermore, $\PRfunc$ has a strong natural number object. That is, for every $f:\N^a \lra \N^b$ and $g:\N^a \times \N^b \times \N \lra \N^b$
there exists a {\it unique} $h=f\sharp g:\N^a \times \N \lra \N^b$ that satisfies the following two commutative diagrams

\be\xymatrix{
\N^a \times \N \ar[rr]^{id \times z} \ar[dd]_{\pi^{a+1}_a} && \N^a \times \N \ar[dd]^h\\
\\
\N^a \ar[rr]_f && \N^b
}
\qquad
\xymatrix{
\N^a \times \N \ar[rr]^{id \times s} \ar[dd]_{\la id, id ,h\ra} && \N^a \times \N \ar[dd]^h\\
\\
\N^a \times \N \times \N^b \ar[rr]_g && \N^b
\label{natnumobj}}\ee

This category of primitive recursive functions was studied extensively by many people including \cite{Burroni,Maietti,mft,Roman,Yano}.
It is known to be the initial object in the 2-category of categories, with products and strict natural number objects.
Other categories in that 2-category will be primitive recursive functions with oracles. One can think of 
oracles as functions put on the leaves of the trees besides the initial functions.

There is a surjective graph morphism $Q:\PRdesc \lra \PRfunc$ that takes $\N^n$ to $\N^n$, i.e., is identity on
objects and $Q$ takes descriptions of primitive recursive functions in $\PRdesc$ to the functions they describe
in $\PRfunc$. Since every primitive recursive function has a --- in fact infinitely many --- primitive recursive description, $Q$ is
surjective on morphisms. Another way to say this is that $\PRfunc$ is a quotient of $\PRdesc$.

Algorithms will be graphs that are ``between'' $\PRdesc$ and $\PRfunc$.

\section{Algorithms}
In the last section we saw the type of structure the set of primitive recursive programs and functions form. In this section we look at the  
types of structures a set of algorithms can have.
\begin{defi}
A {\bf primitive recursive (P.R.) algorithmic universe}, $\PRalg$, is a graph whose objects are the powers of natural
numbers $\N^0, \N^1, \N^2, \ldots, \N^n, \ldots$. We furthermore require that there exist graph morphisms $R$ and $S$ that
are the identity on objects and that make the following diagram of graphs commute:
\be \xymatrix{
\PRdesc\ar@{>>}[dd]_Q \ar@{>>}[rd]^R
\\
& \PRalg \ar@{>>}[ld]^S
\\
\PRfunc.\\
\label{algunidef}}.\ee

The image of the initial functions under $R$ will be distinguished objects in $\PRalg$:
$z:\N^0 \lra \N$, 
$s:\N \lra \N$,  and
for all $n \in \N$ and for all $1\le i \le n$ there
are projection maps $\pi^n_i:\N^n \lra \N$.
\end{defi}

In addition, a P.R. algorithmic universe {\it might} have the following operations: (Warning: even when they exist, these are not necessarily functors because we are not dealing
with categories.)
\begin{itemize}
\item Composition: For $f:\N^a \lra \N^b$ and $g:\N^b \lra \N^c$, there is a $$(g \circ f):\N^a \lra \N^c.$$
\item Recursion: For $f:\N^a \lra \N^b$ and $g:\N^{a+b+1} \lra \N^b$, there is a $$(f \sharp g):\N^{a+1} \lra \N^b $$
\item Bracket: For $f:\N^a \lra \N^b$ and $g:\N^a \lra \N^c$, there is a $$\la f, g \ra:\N^a \lra \N^{b+c}$$
\end{itemize}
These operations are well defined for programs but need not be well defined for equivalence classes of programs. There was
never an insistence that our equivalence relations be congruences (i.e. respect the operations). We study when these operations
exist at the end of the section.

Notice that although the $Q$ graph morphism preserves the composition, bracket and recursion operators, we do not insist that
$R$ and $S$ preserve them. We will see that this is too strict of a requirement.

\begin{defi}
Let $\PRalg$ be a P.R. algorithmic universe. A {\bf P.R. quotient algorithmic universe} is a P.R.
algorithmic universe $\PRalg'$ and an identity on objects, surjection on edges graph map $T$ that
makes all of the following triangles commute
\be \xymatrix{
\PRdesc\ar@{>>}[dd]_Q \ar@{>>}[rr]^R \ar@{>>}[ddrr]^(.3){R'}&&\PRalg \ar@{>>}[dd]^{T}\ar@{>>}[ddll]_(.3){S}
\\
&
\\
\PRfunc&& \PRalg'. \ar@{>>}[ll]^{S'}\label{quotentalg}
}\ee
\end{defi}

\vspace{.2in}

Examples of P.R. algorithmic universe abound:

\begin{examp}
$\PRdesc$ is the primary trivial example. In fact, all our examples will be quotients of this algorithmic universes. Here $R=id$ and $S=Q$.
\end{examp}

\vspace{.2 in}

\begin{examp}
$\PRfunc$ is another trivial example of an algorithmic universe. Here $R=Q$ and $S=id$.
\end{examp}

\vspace{.2 in}

\begin{examp}
$\PRalgC$ is a quotient of $\PRdesc$. This is constructed by adding the following relation:

For any three composable maps $f$, $g$ and $h$, we have
\be h \circ (g \circ f) \sim (h \circ g) \circ f.\label{compisassoc}\ee
In terms of trees, we say that the following trees are equivalent:

\vspace{.2in}

\qtreecenterfalse

\Tree [ [ {$f:\N^a\raa\N^b$} {$g:\N^b\raa\N^c$} ] .{$g\circ f:\N^a\raa\N^c$ \\ \c} {$h:\N^c \raa\N^d$}
].{$h\circ (g\circ f):\N^a\raa\N^d$\\ \c}
$\sim$
\Tree [ {$f:\N^a\raa\N^b$} [ {$g:\N^b\raa\N^c$} {$h:\N^c \raa\N^d$} ].{$h\circ g:\N^b\raa\N^d$\\ \c}
].{$(h\circ g) \circ f:\N^a\raa\N^d$\\ \c}

\vspace{.3in}

It is obvious that if there is a well-defined composition map in $\PRalgC$ it is associative.
\end{examp}

\vspace{.2 in}

\begin{examp}
$\PRalgI$ is also a quotient of $\PRdesc$ that is constructed by adding in the relations that say that the projections $\pi^I$s
act like identity maps.
That means for any $f:\N^a \raa\N^b$,
we have
\be f \circ \pi^a_a \sim f \sim \pi^b_b \circ f.\label{composident}\ee

In terms of trees:
\vspace{.2in}

\qtreecenterfalse

\Tree [{$\pi^{a}_{a}:\N^a \raa\N^a$} {$f:\N^a\raa\N^b$} ].{$f \circ \pi^{a}_{a} :\N^a \raa\N^b$\\ \c}
$\sim$
 {$\quad f:\N^a\raa\N^b \quad$}
$\sim$
\Tree[ {$f:\N^a\raa\N^b$} {$\pi^{b}_{b}:\N^b \raa\N^b$} ].{$\pi^{b}_{b} \circ f:\N^a \raa\N^b $\\ \c}
\vspace{.2in}
\end{examp}

The composition map in $\PRalgI$ has a unit.

\vspace{.2 in}

\begin{examp}
$\PRalgCat$ is $\PRdesc$ with both relations (\ref{compisassoc}) and (\ref{composident}).
Notice that this ensures that $\PRalgCat$ is more than a graph and is, in fact, a full fledged category.
\end{examp}

\vspace{.2 in}
\begin{examp}
$\PRalgCatX$ is a quotient of $\PRalgCat$ which has a well-defined bracket/product function. We add the following relations to $\PRalgCat$:
\begin{itemize}
\item The bracket is associative. For any three maps $f,g,$ and $h$ with the same domain, we have
$$ \la \la f,g \ra, h \ra \sim\la f,\la g, h \ra \ra $$
In terms of trees, this amounts to

\Tree [ [ {$f:\N^a\raa\N^b$} {$g:\N^a\raa\N^c$} ] .{$\la f,g \ra :\N^a\raa\N^b\times \N^c$ \\ \b} {$h:\N^a \raa\N^d$}
].{$\la \la f,g\ra h\ra:\N^a\raa\N^b\times \N^c\times \N^d$\\ \b}
$\sim$
\Tree [ {$f:\N^a\raa\N^b$} [ {$g:\N^a\raa\N^c$} {$h:\N^a \raa\N^d$} ].{$\la g,h \ra:\N^b\raa\N^c \times \N^d$\\ \b}
].{$\la f, \la g,h \ra \ra:\N^a\raa\N^b\times \N^c\times \N^d$\\ \b}

\vspace{.2in}

\item Composition distributes over the bracket on the right.
For $g:\N^a \raa\N^b$,  $f_1:\N^b \raa\N^c$ and $f_2:\N^b \raa\N^d,$ we have
\be \la f_1,f_2 \ra \circ g \sim \la f_1 \circ g, f_2 \circ g \ra.\label{distrib}\ee

In terms of trees, this amounts to saying that these trees are equivalent:

\Tree [ {$g:\N^a\raa\N^b$} [ {$f_1:\N^b\raa\N^c$} {$f_2:\N^b\raa\N^d$}
 ].{$\la f_1, f_2 \ra:\N^b\raa\N^c \times \N^d$\\ \b} ].{$\la f_1, f_2 \ra \circ g:\N^a\raa\N^c \times \N^d$\\ \c}
\Tree [ [{$g:\N^a\raa\N^b$} {$f_1:\N^b\raa\N^c$} ].{$f_1 \circ g:\N^a\raa\N^c$\\ \c}
[{$g:\N^a\raa\N^b$} {$f_2:\N^b\raa\N^d$} ].{$f_2 \circ g:\N^a\raa\N^d$\\ \c}
].{$\la f_1\circ g, f_2\circ g \ra:\N^a\raa\N^c \times \N^d$\\ \b}

\vspace{.2in}
\
\item The bracket is almost commutative.
For any two maps $f$ and $g$ with the same domain,
$$\la f,g \ra \sim tw \circ \la g, f \ra.$$
In terms of trees, this amounts to

\Tree [{$f:\N^a\raa\N^b$} {$g:\N^a\raa\N^c$} ] .{$\la f,g \ra :\N^a\raa\N^b\times \N^c$ \\ \b}
$\sim$
\Tree [ [ {$g:\N^a\raa\N^c$} {$f:\N^a\raa\N^b$} ] .{$\la g,f \ra :\N^a\raa\N^c\times \N^b$ \\ \b} {$tw:\N^c \times \N^b\raa\N^b \times\N^c$}
].{$tw \circ \la g,f \ra: \N^a \raa\N^b \times \N^c$ \\ \c}

\vspace{.2in}

\item Twist is idempotent.
$$tw_{\N^a,\N^b} \circ tw_{\N^a,\N^b} \sim id = \pi^{a + b}_{a +b}:\N^a \times \N^b \raa\N^a \times \N^b.$$

\vspace{.2in}
\item Twist is coherent. That is, the twist maps of three elements behave with respect to themselves.
$$(tw_{\N^b,\N^c}\times id) \circ (id \times tw_{\N^a, \N^c})\circ (tw_{\N^a,\N^b}\times id)
\sim
(id \times tw_{\N^a,\N^b})\circ (tw_{\N^a,\N^c}\times id)\circ (id \times tw_{\N^b,\N^c}).$$
This is called the hexagon law or the third Reidermeister move. Given the idempotence and hexagon laws,
it is a theorem
that there is a unique twist map made of smaller twist maps between any two products of elements (\cite{MacLane} Section
XI.4). The induced product map will be coherent.
\end{itemize}
\end{examp}

\vspace{.2 in}
\begin{examp}
$\PRalgCatN$ is a category with a natural number object. It is $\PRalgCat$ with the following relations:
\begin{itemize}
\item Left square of Diagram (\ref{natnumobj}). $(f\sharp g) \circ(id \times z) \quad \sim \quad (f \circ \pi^{a+1}_a).$
\item Right square of Diagram (\ref{natnumobj}). $ (f\sharp g) \circ(id \times s) \quad \sim \quad (f \circ \la id, id, (f \sharp g) \ra .$
\item Natural number object and identity. If $g=\pi^{a+b+1}_b : \N^a \times \N^b \times \N \lra \N^b$ then
$$(f\sharp \pi^{a+b+1}_b) \quad \sim \quad (f \circ \pi^{a+1}_a).$$
\item Natural number object and composition. This is explained in Section 3.5 of \cite{Yano}.
$$g_1 \circt (f \sharp (g_2 \circt g_1)) \quad \sim  \quad (g_1 \circt f)\sharp (g_1 \circt g_2).$$
\end{itemize}
\end{examp}
\vspace{.2 in}

\begin{examp}
$\PRalgCatXN$ is a category that has both a product and natural number object. It can be constructed by
adding to $\PRalgCat$ all the relations of $\PRalgCatX$ and $\PRalgCatN$ as well as the following relations:
\begin{itemize}
\item Natural number object and bracket. This is explained in Section 3.4 of \cite{Yano}
$$\la f_1,f_2 \ra \sharp (g_1 \boxtimes g_2) \quad \sim \quad  \la f_1 \sharp g_1, f_2 \sharp g_2 \ra.$$
\end{itemize}
\end{examp}

Putting all these examples together, we have the following diagram of P.R. algorithmic universes.

$$\xymatrix{
&&\PRdesc \ar@{>>}[ld] \ar@{>>}[dr]
\\
& \PRalgC \ar@{>>}[dr]&&\PRalgI \ar@{>>}[dl]
\\
&& \PRalgCat \ar@{>>}[ld] \ar@{>>}[dr]
\\
& \PRalgCatX \ar@{>>}[dr]&&\PRalgCatN \ar@{>>}[dl]
\\
&& \PRalgCatXN \ar@{>>}[d]
\\
&& \PRfunc\\
}$$

There is no reason to think that this is a complete list. One can come up with infinitely many more examples of algorithmic universes. We can take other permutations and combinations of the relations given here as well as new ones. Every appropriate equivalence relation will give a different algorithmic universe. 

In \cite{Yano}, we mentioned other relations which deal with the relationship between the operations and the initial functions. We
do not mention those relations here because our central focus is the existence of well defined operations.

\vspace{.2in}

A word about decidability. The question is, for a given  P.R. algorithmic universe determine whether or not
two programs in $\PRdesc$ are in the same equivalence class of that algorithmic universe.
\begin{itemize}
\item This is very easy in
the algorithmic universe $\PRdesc$ since every equivalence class has only one
element. Two descriptions are in the same equivalence relation iff they are {\it exactly} the same.
\item The extreme opposite is
in $\PRfunc$. By a theorem similar to Rice's theorem, there is no way to tell when two different programs/descriptions are the
same primitive recursive function. So $\PRfunc$ is not decidable.
\item In between $\PRdesc$ and $\PRfunc$ things get a little hairy. This is the boundary between syntax and semantics. 
Consider $\PRalgC$, i.e., the graph with associative composition. This is decidable. All one has to do is change all the contiguous sequences
of compositions to associate on the left. Do this for both descriptions and then see if the two modified programs are the same.
\item One can perform a similar trick for $\PRalgI$. Simply eliminate all the identities and see if the two modified programs are the same.
\item For $\PRalgCat$ one can combine the tricks from $\PRalgC$ and $\PRalgI$ to show that it is also decidable.
\item $\PRalgCatX$ is also decidable because of the coherence of the product. Once again, any contiguous sequences of products
can be associated to the left. Also, equivalence relation (\ref{distrib}) insures the naturality of the product so that products and compositions can
``slide across'' each other. Again, each description can be put into a canonical form and then see if the modified programs are the same.
\item However we loose decidability when it comes to structures with natural number objects. See the important paper by Okada and Scott \cite{Okada}. It seems that this implies that $\PRalgCatN$ and $\PRalgCatXN$ are undecidable.  One can think of this as
the dividing line between the decidable, syntactical structure of $\PRdesc$ and the undecidable, semantical structure of $\PRfunc$.
\end{itemize}

\section{Galois Theory}

An automorphism $\phi$ of $\PRdesc$ is a graph isomorphism that is the identity on the vertices (i.e., $\phi(\N^n)=\N^n)$. For every $a,b \in \N$ $\phi$ basically acts on the edges between 
$\N^a$ and $\N^b$. 
We are interested in automorphisms that preserve functionality. That is, automorphisms $\phi$, such that for all programs $p$,
 we have that $p$ and $\phi(p)$ perform the same function. In terms of Diagram (\ref{GalofAlg}) we demand that $Q(\phi(p))=Q(p)$.
It is not hard to see that the set of all automorphism of $\PRdesc$ that preserve functionality
forms a group. We denote this group as $Aut(\PRdesc/\PRfunc)$. We shall look at subgroups of this group and see its relationship with intermediate fields. Let $GRP$ denote the partial order of subgroups of $Aut(\PRdesc/\PRfunc)$

Let $ALG$ be the partial order of intermediate algorithmic universes. One algorithmic universe, $\PRalg$ is greater than or equal to  $\PRalg'$ if there is a quotient algorithmic  map $\PRalg \twoheadrightarrow \PRalg'$.

We shall construct a Galois connection between $GRP$ and $ALG$. That is, there will be an order reversing map $\Phi: ALG \lra GRP$ and an order reversing $\Psi: GRP \lra ALG$.

In detail, for a given algorithmic universe $\PRalg$, we construct the subgroup 
$$\Phi(\PRalg)\subseteq Aut(\PRdesc/\PRfunc).$$
$\Phi(\PRalg)$ is the set of all automorphisms of $\PRdesc$ that preserve that algorithmic
universe, i.e., automorphisms $\phi$ such that for all programs $p$, we have $p$ and $\phi(p)$ are in
the same equivalence class in $\PRalg$. That is,
$$\Phi(\PRalg) = \{\phi| \forall p \in \PRdesc, \quad [\phi(p)]=[p] \in \PRalg \}.$$
In terms of Diagram (\ref{algunidef}), this means $R(\phi(p))=R(p)$.
In order to see that $\Phi(\PRalg)$ is a subgroup of $Aut(\PRdesc/\PRfunc)$, notice that
if $\phi$ is in $\Phi(\PRalg)$ then we have
$$R\phi=R \qquad \Rightarrow \qquad SR\phi=SR \qquad \Rightarrow \qquad Q\phi=Q$$
which means that $\phi$ is in $Aut(\PRdesc/\PRalg)$.
In general, this subgroup fails to be normal.

If $T:\PRalg \twoheadrightarrow \PRalg'$ is a quotient algorithmic universe as in Diagram (\ref{quotentalg}) then
$$\Phi(\PRalg) \subseteq \Phi(\PRalg')\subseteq Aut(\PRdesc/\PRfunc).$$
This is obvious if you look at $\phi \in \Phi(\PRalg)$ then we have that
$$R\phi = R \qquad \Rightarrow \qquad TR\phi=TR \qquad \Rightarrow \qquad R'\phi=R'$$
which means that $\phi$ is also in $\Phi(\PRalg').$

\vspace{.3in}

The other direction goes as follows. For $H \subseteq Aut(\PRdesc/\PRfunc)$, the graph $\Psi(H)$ is a quotient of $\PRdesc$.
The vertices of $\Psi(H)$
are powers of natural numbers. The edges will be equivalence classes of edges from $\PRdesc$. The equivalence relation $\sim_H$
is defined as
\be p \sim_H p' \qquad \mbox{iff} \qquad \mbox{there exists a } \phi \in H \mbox{ such that }\phi(p)=p' \label{defofequiv}\ee
The fact that $\sim_H$ is an equivalence relation follows from the fact that $H$ is a group. In detail
\begin{itemize}
\item Reflexivity comes from the fact that $id\in H$.
\item Symmetry comes from the fact that if $\phi \in H$ then $\phi^{-1} \in H$.
\item Transitivity comes from the fact that if $\phi \in H$ and $\psi \in H$ then $\phi\psi\in H$.
\end{itemize}

If $H \subseteq H' \subseteq Aut(\PRdesc/\PRfunc)$ then there will be a surjective map
$$\Psi(H) \twoheadrightarrow \Psi(H').$$
The way to see this is to realize that there are more $\phi$ in $H'$ to make different programs
equivalent as described in line (\ref{defofequiv}).

\vspace{.3in}

\begin{teo} The maps $\Phi: ALG \lra GRP$ and $\Psi: GRP \lra ALG$ form a Galois connection. \end{teo}

\noindent {\bf Proof.} We must show that for any $H$ in $GRP$ and any $\PRalg$ in $ALG$ we have 
$$ H \subseteq \Phi (\PRalg ) \mbox{ if and only if } \Psi(H)  \twoheadrightarrow \PRalg.$$
This will be proven with the following sequence of implications.
$$ H \subseteq \Phi (\PRalg )$$
$$ \mbox{ if and only if}$$
$$ \phi \in H  \quad \Rightarrow \quad  \phi \in \Phi (\PRalg )$$
$$ \mbox{ if and only if}$$
$$ \phi \in H \quad \Rightarrow \quad \forall p, [\phi(p)]=[ p] \in \PRalg $$
$$ \mbox{ if and only if}$$
$$\forall p  \quad \phi(p) \sim_H p\quad \Rightarrow \quad \phi(p) \sim_{\PRalg} p$$
$$ \mbox{ if and only if}$$
$$  \Psi(H)  \twoheadrightarrow \PRalg.$$  $\Box$

\vspace{.3in}

Every Galois connection (adjoint functor) induces an isomorphism of sub-partial orders (equivalence of categories.) Here we do not have to look at a sub-partial order of $ALG$ for the following reason:

\begin{teo} For any $\PRalg$ in $ALG$ 
$$(\Psi \circ \Phi)(\PRalg)=\PRalg.$$
\end{teo}

\noindent {\bf Proof.} $\Psi ( \Phi((\PRalg))= \sim_{ \Phi((\PRalg)}$ $$=\sim_{\Phi(\{\phi | \forall p \in \PRdesc, \quad [\phi(p)]=[p] \in \PRalg \})}=\PRalg.\Box$$

In contrast, it is not necessarily the case that for any $H$ in $GRP$, we have 
$$ (\Phi \circ \Psi )(H) = H.$$
We do have that 
$$ (\Phi \circ \Psi )(H)= \Phi(\Psi(H))=\Phi(\sim_{H})=\{\phi| \forall p \quad \phi(p)\sim_H p \} \supseteq H$$
because any $\phi$ in $H$ definitely satisfies that condition. But many other $\phi$ might also satisfy this requirement.  
In general this requirement is not satisfied. $H$ might generate a transitive action. In that case $\Phi(\Psi(H))$ will be all automorphisms.


A subgroup $H$ whose induced action does not extend beyond $H$ will be important:
\begin{defi}
A subgroup $H$ of $Aut(\PRdesc/\PRfunc)$ is called ``restricted'' if $ (\Phi \circ \Psi )(H) = H.$
\end{defi}

We can sum up with the following statement:

\begin{teo}[Fundamental theorem of Galois theory] The lattice of restricted subgroups of $Aut(\PRdesc/\PRfunc)$ is isomorphic to the dual
lattice of algorithmic universes between $\PRdesc$ and $\PRfunc$.
\end{teo}

Notice that the algorithmic universes that we dealt with in this theorem does not necessarily have well-defined 
extra structure/operations.
We discussed the equivalence relations of $\PRdesc$ and did not discuss congruences of
$\PRdesc$. Without the congruence, the operations of composition, bracket and recursion might not be well-defined for the equivalence classes. This is very similar to classical Galois theory where we discuss a single weak structure (fields) and discuss all intermediate objects as fields even though they might have more structure. So too, here we stick to a weak structure. 

However we can go further. Our definition of algorithmic universes is not carved in stone. One can go on and define, say,  a {\it composable} algorithmic universe. This is an algorithmic universe with a well-defined composition function. Then we can make the fundamental theorem of Galois theory for composable algorithmic universes by looking at automorphisms of $\PRdesc$ that preserve the composition operations. That is, automorphisms $\phi$ such that 
for all programs $p$ and $p'$ we have that $\phi(p\circ p')=\phi(p) \circ \phi(p').$
Such automorphisms also form a group and one can look at subgroups as we did in the main
theorem. On the algorithmic universe side, we will have to look at equivalence relations that are 
congruences. That is, $\sim$ such that if $p\sim p'$ and $p'' \sim p'''$ then
$p \circ p''  \sim  p' \circ p'''.$
Such an analogous theorem can easily be proved. 

Similarly, one can define {\it recursive} algorithmic universes and {\it bracket} algorithmic universes.  
One can still go further and ask that an algorithmic universe has two well-defined operations. 
In that case the automorphism will have to preserve two operations. If $H$ is a group of automorphisms, then we can denote the subgroup of automorphisms that preserve composition as 
$H_C$. Furthermore, the subgroup that preserves composition and recursion will be denoted as $H_{CR}$, etc. The subgroups fit into the following lattice.

\be \xymatrix{
&&H
\\
\\
H_C\ar@{^{(}->}[uurr]&&H_R\ar@{^{(}->}[uu]&&H_B\ar@{^{(}->}[uull]
\\
\\
H_{CR}\ar@{^{(}->}[uu]\ar@{^{(}->}[uurr]|\hole &&H_{CB}\ar@{^{(}->}[uull]\ar@{^{(}->}[uurr]
&&H_{RB}\ar@{^{(}->}[uull]|\hole \ar@{^{(}->}[uu]
\\
\\
&&H_{CRB}=\{e\}.\ar@{^{(}->}[uull]\ar@{^{(}->}[uurr]\ar@{^{(}->}[uu]
}\ee

It is important to realize
that it is uninteresting to require that the algorithmic universe have all three operations. 
The only automorphism that preserves all the operations is the identity automorphism on $\PRdesc$. One can see this by remembering that the automorphisms preserve all the initial functions and if we 
ask them to preserve all the operations, then it must be the identity automorphism. This is similar to looking at an automorphism of a group that preserves the generators and the group operation. That is not a very interesting automorphism. 

One can ask of the automorphisms to preserve all three operations but not preserve the initial operations. Similarly, when discussing oracle computation, one can ask the automorphisms to preserve all three operations and the initial functions, but not the oracle functions. All these suggestions open up new vistas of study.   

\section{Homotopy Theory}
The setup of the structures presented here calls for an analysis from the homotopy perspective. We seem to have a covering space and are looking at automorphims of that covering space. Doing homotopy theory from this point of view makes it very easy to generalize to many other areas in mathematics and computer science. This way of doing homotopy theory is sometimes called Grothendieck's Galois theory. We gained much from \cite{Baez} and \cite{Dubuc}. 

First a short review of classical homotopy theory. For $X$ a ``nice'' connected space and $P:C \lra X$ a universal covering space, it is a fact that $P$ induces an isomorphism 
\be Aut(C) \qquad \cong \qquad \pi_1(X)\ee 
where $Aut(C)$ is the set of automorphisms (homeomorphisms) of $C$ that respect $P$ (automorphisms $f: X \lra X$ such that $P \circ f= P$.) Such automorphisms are called ``deck transformations.'' $\pi_1(X)$ is the fundamental group of $X$. This result can be extended by dropping the assumption that $X$ is connected. We then have  
\be Aut(C) \qquad \cong \qquad \Pi_1(X)\ee
where $\Pi_1(X)$ is the fundamental groupoid of $X$. Another way of generalizing this result is to consider $C$ which is not necessarily a universal covering space. In other words, consider $C$ where $\pi_1(C)$ is not necessarily trivial. The theorem then says that 
\be Aut(C) \qquad \cong \qquad \pi_1(X)/P_*(\pi_1(C)).\ee
That is, we look at quotients of $\pi_1(X)$ by the image of the fundamental groups of $C$.

Our functor $Q: \PRdesc \lra \PRfunc$ seems to have a feel of a covering space. The functor is bijective on objects and is surjective on morphisms. Also, for every primitive recursive function $f:\N^m \lra \N^n$ the set 
$$P^{-1}(f:\N^m \lra \N^n)$$
has a countably infinite number of programs/descriptions.

Our goal will be an isomorphism of the form 
\be Aut(\PRdesc/\PRfunc) \qquad \cong \qquad \pi_1(\PRfunc)\ee  
The right side of the purported isomorphism is very easy to describe. $\pi_1(\PRfunc)$ is the group of invertible primitive recursive functions from $\N$ to $\N$. Note that because of  primitive recursive isomorphisms of the form $\N \lra \N^k$ for all $k$ (G\"{o}del numbering functions), the elements of this group can be rather sophisticated.

One should not be perturbed by the fact that we are looking at the reversible primitive recursion functions based at $\N$ as opposed to $\N^k$ for an arbitrary $k$. We can also look at the fundamental group based at $\N^k$ and denote these two fundamental groups as $\pi_1(\PRfunc, \N)$ and $\pi_1(\PRfunc, \N^k)$. By a usual trick of classical homotopy theory, these two groups are isomorphic as follows. Let $\alpha:\N^k \lra \N$ be a primitive recursion isomorphic function. Let $f: \N \lra \N$ be an element of $\pi_1(\PRfunc, \N)$. We then have the isomorphism from $\pi_1(\PRfunc, \N)$ to $\pi_1(\PRfunc, \N^k)$:
$$f: \N \lra \N \qquad \mapsto \qquad \alpha^{-1} \circ f \circ \alpha: \N^k \lra \N \lra \N \lra \N^k.$$
Since all these groups are isomorphic, we ignore the base point and call the group $\pi_1(\PRfunc)$.

As far as I can find, this group has not been studied in the literature. There are many questions to ask. Are there generators to this group? What properties does this group have?(An anonymous reviewer gave a simple proof that it is not finitely generated.)  What is the relationship of the groups of all recursive isomorphisms from $\N$ to $\N$ and the primitive recursive isomorphisms? What is the relationship between invertible primitive recursive functions (that is, a primitive recursive function that is an isomorphism) and all primitive recursive functions? Can every primitive recursive function be mimicked in some way (like Bennett’s result about reversible computation) by a reversible/invertible primitive recursive function?  In essence, the group $\pi_1(\PRfunc)$ is the upper left of the following commutative square of monoides:
$$\xymatrix{\txt{group of\\ invertible primitive \\ recursive functions} \ar@{^{(}->}[rr]\ar@{^{(}->}[dd]&& \txt{group of\\ invertible  \\ recursive functions}\ar@{^{(}->}[dd]
\\
\\
\txt{monoid of\\ all primitive \\ recursive functions}\ar@{^{(}->}[rr] && \txt{monoid of\\  recursive functions}}$$
Essentially we are asking if any of these inclusion maps have some type of retract.

Unfortunately the map $Q: \PRdesc \lra \PRfunc$ fails to be a real covering map because it does not have the ``unique lifting property.'' In topology, a covering map $P:C \lra X$ has the following lifting property: for any path $f:I \lra X$ from $f(0)=x_0$ to $f(1)=x_1$ and for any $c_0 \in C$ such that $P(c_0)=x_0$ there is a {\em unique} $\hat{f}:I \lra C$ such that $P(\hat{f})=f$. In English this says that for any path in $X$ and any starting point in $C$, there is a unique path in $C$ that maps onto the path in $X$. In our context, such a unique lifting would mean that for every primitive recursive function made out of a sequence of functions, if you choose one program/description to start the function, then the rest of the programs/descriptions would all be forced. There is, at the moment, no reason for this to be true. The problem is that above  
every function, there is only a set of programs/descriptions. This set does not have any more structure. 

However, all hope is not lost. Rather than look at $\PRdesc$ as simply a graph, look at it as graph enriched in groupoids. That is, between every two edges there is a possibility for there to be isomorphisms corresponding to whether or not the two programs are essentially the same. This is almost a bicategory or 2-groupoid but the bottom level is not a category but a graph. (I am grateful to Robert Par\'e for this suggestion.) If we were able to formalize this, then we would have that 
$P^{-1}(f:\N^m \lra \N^n)$ would be a connected groupoid and we would have some type of unique lifting property. At the moment I am not sure how to do this. Another advantage of doing this is that $$Aut(\PRdesc/\PRfunc)$$ would not be as big and as  unmanageable as it is now. The automorphisms would have to respect the higher dimensional cells. Much work remains. 

\section{Future Directions}

\noindent {\bf Extend to all computable functions.} The first project worth doing is to extend this work to all computable functions from primitive recursive functions. One need only add in the minimization operator and look at its relations with the other operations.
The study of such programs from our point is an ongoing project in \cite{ManinYano}. However the Galois theory perspective
is a little bit complicated because of the necessity to consider partial operations. A careful study of \cite{Rosenberg} will,
no doubt, be helpful.

\vspace{.2 in}

\noindent {\bf Continuing with Galois Theory} There are many other classical Galois theory theorems that need to
be proved for our context. We need
the Zassenhaus lemma,
the Schreier refinement theorem, and culminating in the
Jordan-H\"older theorem. In the context of algorithms this would be some statement about decomposing a category of algorithms regardless of the order in which the equivalence relations are given. We might also attempt a form of the
Krull-Schmidt theorem.

\vspace{.2 in}

\noindent{\bf Impossibility results.} The most interesting part of Galois theory is that it shows that there
are certain contingencies that are impossible or not ``solvable.'' What would the analogue for algorithms be?

\vspace{.2 in}

\noindent{\bf Calculate some groups.} It is not interesting just knowing that there are automorphism groups.
It would be nice to actually give generators and relations for some of these groups. This will also give us a firmer
grip for any impossibility results.

\vspace{.2 in}

\noindent{\bf Universal Algebra of Algorithms.} In this paper we
stressed looking at quotients of the structure of all programs.
However there are many other aspects of the algorithms that we can
look at from the universal algebraic perspective. \underline{Subalgebras:}
We considered all primitive recursive programs. But there are
subclasses of programs that are of interest. We can for example
restrict the number of recursions in our programs and get to
subclasses like the Grzegorczyk's hierarchy. How does the subgroup lattice survive with this stratification? Other subclasses of primitive
recursive functions such as polynomial functions and EXPTIME
functions can also be studied. \underline {Superalgebras:} We can also look
at larger classes of algorithms. As stated above, we can consider
all computable functions by simply adding in a minimization
operator. Also, oracle computation can be dealt with by looking at
trees of descriptions that in addition to initial functions permit
arbitrary functions on their leaves. Again we ask similar questions
about the structure of the lattice of automorphisms and the related lattice of intermediate algorithms.
\underline{Homomorphisms:} What would correspond
to a homomorphism between classes of computable algorithms?
Compilers. They input programs and output programs. This opens up a
whole new can of worms. What does it mean for a compiler to preserve
algorithms? When are two compilers similar? What properties should a
compiler preserve? How are the lattices of subgroups and intermediate algorithms preserved under homomorphisms/compilers?
There is obviously much work to be done.

\end{document}